.
.
.
\font\sets=msbm10.
\font\stampatello=cmcsc10.
\font\symbols=msam10.

\def\1{{\bf 1}}

\def\starsum{\mathop{\enspace{\sum}^{\ast}}}

\def\square{\hbox{\vrule\vbox{\hrule\phantom{s}\hrule}\vrule}}
\def\defineq{\buildrel{def}\over{=}}
\def\defin{\buildrel{def}\over{\Longleftrightarrow}}
\def\doublesum{\mathop{\sum\sum}}
\def\doublestarsum{\mathop{\starsum\starsum}}

\def\supporto{{\rm supp}\,}
\def\C{\hbox{\sets C}}
\def\N{\hbox{\sets N}}

\def\R{\hbox{\sets R}}
\def\Z{\hbox{\sets Z}}

\def\EssBdd{\hbox{\symbols n}\,}
\def\modSel{{\widetilde{J}}}
\def\modchi{{\widetilde{\chi}}}
\def\modcjq{{\widetilde{c}_{j,q}}}
\def\modFh{{\widetilde{F}}_h\!\!}

\par
\centerline{\bf ON THE MODIFIED SELBERG INTEGRAL}
\bigskip
\centerline{Giovanni Coppola}
\bigskip
{
\font\eightrm=cmr8
\eightrm {
\par
{\bf Abstract.} We give a kind of \lq \lq approximate majorant principle\rq \rq \thinspace result for the \lq \lq modified Selberg integral\rq \rq, say $\modSel_f(N,h)$, of essentially bounded $f:\N \rightarrow \R$ (i.e., bounded by arbitrary small powers); i.e., we get an upper bound, in terms of the modified Selberg integral of a related function $F$ (with $|f\ast \mu|\ll F\ast \mu$, in the supports intersection), getting a \lq \lq square-root cancellation\rq \rq \thinspace  for the error-terms. Here $\modSel_f(N,h)$ is the mean-square (in $N<x\le 2N$) of the \lq \lq averaged short sum\rq \rq \thinspace of, say, $f:=g\ast \1$, minus its expected value; i.e., ${1\over h}\sum_{m\le h}\sum_{0\le |n-x|<m}f(n)-M_f(x,h)$, with expected value $M_f(x,h)$ (say, $\approx h\sum_{d\le x}g(d)/d$); so, this mean-square weights, on average, the $f-$values in (almost all, i.e. all, but $o(N)$ possible exceptions) the short intervals $[x-h,x+h]$, with mild restrictions on \thinspace $h$ \thinspace (say, \thinspace $h\to \infty$ \thinspace and \thinspace $h=o(N)$, when \thinspace $N\to \infty$).
}
\footnote{}{\par \noindent {\it Mathematics Subject Classification} $(2010) : 11N37, 11N25.$}
}
\bigskip
\par
\noindent \centerline{\bf 1. Introduction and statement of the results.}
\smallskip
\par
\noindent
We give upper bounds for the {\stampatello modified Selberg integral} (for the Selberg integral [C-S], [C1], [C2], [C3])
$$
\modSel_f(N,h)\defineq \sum_{x\sim N} \Big| \sum_{|n-x|\le h}\Big(1-{{|n-x|}\over h}\Big)f(n)-M_f(x,h)\Big|^2 
= \sum_{x\sim N} \Big| {1\over h}\sum_{m\le h} \sum_{0\le |n-x|<m}f(n)-M_f(x,h)\Big|^2, 
$$
\par
\noindent
(now on $x\sim N$ is $N<x\le 2N$ in sums), where the {\stampatello mean-value} (averaged sum's expected value) is ($\forall \varepsilon>0$) 
$$
M_f(x,h)\defineq h\sum_{d\le x+h}{{g(d)}\over d}
=h\sum_{d\le x}{{g(d)}\over d}+{\cal O}\Big( {h\over x}\sum_{x<d\le x+h}|g(d)|\Big) 
=h\sum_{d\le x}{{g(d)}\over d}+{\cal O}_{\varepsilon}\Big( {{h^2x^{\varepsilon}}\over x}\Big), 
$$
\par
\noindent
with $g:=f\ast \mu$ (see [T]), here for the class of {\stampatello essentially bounded} arithmetic real functions $f$. We use \lq \lq {\stampatello essentially}\rq \rq \thinspace to leave (as they're negligible) arbitrarily small powers of $N$. With Vinogradov notation [D] 
$$
f \enspace \hbox{\rm is } \hbox{\stampatello essentially bounded } (\hbox{\rm abbrev. } f\EssBdd 1) \enspace \defin \enspace \forall \varepsilon>0 \enspace f(n)\ll_{\varepsilon} n^{\varepsilon}
$$
\par
\noindent
(here, $\forall n\le 2N+h$: we \lq \lq don't see\rq \rq \thinspace $f$ any further; hence, $f(n)\ll_{\varepsilon} N^{\varepsilon}$), while $G$ essentially bounds $F$ when 
$$
F(N,h)\EssBdd G(N,h) \enspace \defin \enspace \forall \varepsilon>0 \enspace |F(N,h)|\ll_{\varepsilon} N^{\varepsilon}G(N,h).
$$
\par
\noindent
As an application, ($[x]$ is the {\stampatello integer part} of $x\in \R$) \thinspace $f\EssBdd 1$ \thinspace ($\Leftrightarrow g\EssBdd 1$, by M\"obius inversion [T]) gives 
$$
\int_{N}^{2N}\Big| \sum_{|n-x|\le h}\Big(1-{{|n-x|}\over h}\Big)f(n)-M_f(x,h)\Big|^2dx \EssBdd \int_{N}^{2N}\Big| \sum_{|n-[x]|\le h}\Big(1-{{|n-[x]|}\over h}\Big)f(n)-M_f([x],h)\Big|^2dx + N 
$$
$$
\EssBdd \sum_{N\le x<2N}\Big| \sum_{|n-x|\le h}\Big(1-{{|n-x|}\over h}\Big)f(n)-M_f(x,h)\Big|^2 + N
\EssBdd \modSel_f(N,h) + h^2 + N. 
$$
\par
\noindent
Hence, leaving $\EssBdd N+h^2$, this integral (continuous mean-square) bound comes from the one for \thinspace $\modSel_f(N,h)$. 

\medskip

\par
Our main result is the following. (We call the $G\ast \1$ in the following a {\stampatello Wintner majorant} of $g\ast \1$.)
\smallskip
\par
\noindent {\stampatello Theorem.} {\it Let } $N,h,Q\in \N$, {\it with} {\stampatello even } $h\to \infty$, $h=o(N)$ {\it when } $N\to \infty$ {\it and } $Q\le N+h$. {\it Let } $g:\N \rightarrow \R$ {\it have } $\supporto(g)\subset [1,Q]$. {\it Then} 
$$
|g|\ll G \, \EssBdd 1 \enspace (\hbox{\stampatello in \thinspace the \thinspace supports \thinspace intersection}) \enspace \Rightarrow \enspace \modSel_{g\ast \1}(N,h)\EssBdd \modSel_{G\ast \1}(N,h)+Nh. 
$$

\medskip
\par				
The same arguments in $\S3$ prove the 
\smallskip
\par
\noindent {\stampatello Proposition.} {\it Let } $N,h\in \N$, {\it with} {\stampatello even } $h\to \infty$ {\it and } $h=o(N)$ {\it when } $N\to \infty$. {\it  Assume } $Q:\R \rightarrow \R$ {\it is a strictly increasing function, with } $1\le Q(x)\le x$, $\forall x\ge 1$. {\it Let} \thinspace $g$ {\it be} {\stampatello real} {\stampatello and} ({\it independent of $N,h$ and $x$}) {\it with support} {\stampatello restricted to}, {\it inside $x-$mean-square, the interval} \thinspace $[1,Q(x+h)]$, $\forall x\sim N$. {\it Then} 
$$
|g|\ll G \, \EssBdd 1 \enspace (\hbox{\stampatello in \thinspace the \thinspace supports \thinspace intersection}) \enspace \Rightarrow \enspace \modSel_{g\ast \1}(N,h)\EssBdd \modSel_{G\ast \1}(N,h)+Nh. 
$$

\par
\noindent 
We start proving the Theorem. (The proof will be completed in $\S3$, where we'll prove the Proposition, too.)

\smallskip

\par
\noindent {\stampatello proof.}$\!$ Assume now on $Q\le N+h$, $g:\N \rightarrow \R$, with $\supporto(g)\subset [1,Q]$, $g\ast \1:=f\EssBdd 1$, $h\to \infty$ and $h=o(N)$. 
$$
\modSel_f(N,h)=\sum_{x\sim N}\Big| \sum_{q\le Q}g(q)\modchi_q(x)\Big|^2,
$$
\par
where we define, this time (compare [C-S], esp.), $\forall q\in \N$, (even if we need it $\forall q\le Q$)
$$
\modchi_q(x)\defineq \sum_{{|n-x|\le h}\atop {n\equiv 0(\!\!\bmod q)}}\Big(1-{{|n-x|}\over h}\Big)-{h\over q}=\sum_{{\ell|q}\atop {\ell>1}}{{\ell}\over q}\starsum_{j\le {{\ell}\over 2}}\widetilde{c}_{j,\ell}\cos {{2\pi xj}\over {\ell}}, \enspace \hbox{\stampatello where}
$$
\par
the $\underline{\hbox{\stampatello Fourier coefficients}}$ {\stampatello are positive} (better, non-negative), say from Fej\'er's kernel, 
$$
\modcjq:={1\over q}\Big({2\over h}\,{{\sin^2 \pi jh/q}\over {\sin^2 \pi j/q}}\Big):={1\over q}\modFh\,\Big({j\over q}\Big)\ge 0 \quad \forall j\le {q\over 2}
$$
\par
and (use Parseval identity [C-S]), writing henceforth \enspace $\starsum$ \enspace to {\stampatello sum over reduced residue classes}, 
$$
\starsum_{j\le q}|\modcjq|^2 \ll \sum_{j\le q}|\modcjq|^2 \ll \left\Vert {h\over q}\right\Vert \ll \min\left(1,{h\over q}\right).
$$
\par
Here, as usual, $\Vert \alpha \Vert:=\min_{n\in \Z}|\alpha-n|$ \thinspace is the {\stampatello distance to} the {\stampatello integers}, $\forall \alpha \in \R$. (In fact, 
$$
\modchi_q(x)={1\over q}\sum_{j\not \equiv 0(\!\!\bmod q)}\sum_{|s|\le h}(1-|s|/h)e_q(js)e_q(xj), 
$$
\par
from {\stampatello orthogonality of additive characters} [V]; from $h$ {\stampatello even}, use Fej\'er kernel summation, 
$$
\modchi_q(x)={1\over q}\sum_{0<|j|\le q/2}\Big({1\over h}\,{{\sin^2 \pi jh/q}\over {\sin^2 \pi j/q}}\Big)\cos {{2\pi xj}\over q}
={1\over q}\sum_{j\le q/2}\Big({2\over h}\,{{\sin^2 \pi jh/q}\over {\sin^2 \pi j/q}}\Big)\cos {{2\pi xj}\over q} = 
$$
$$
=\sum_{{d|q}\atop {d<q}}\sum_{{j\le q/2}\atop {(j,q)=d}}\modcjq \cos {{2\pi xj}\over q}
=\sum_{{d|q}\atop {d<q}}{1\over d}\sum_{{j'\le q/(2d)}\atop {(j',(q/d))=1}}\widetilde{c}_{j',q/d}\cos {{2\pi xj'}\over {q/d}}
=\sum_{{\ell|q}\atop {\ell>1}}{{\ell}\over q}\sum_{{j\le \ell/2}\atop {(j,\ell)=1}}\widetilde{c}_{j,\ell}\cos {{2\pi xj}\over {\ell}}, 
$$
\par
where, as in [C-S], we use that \thinspace $\widetilde{c}_{dj',dq'}=\widetilde{c}_{j',q'}/d$, $\forall d,j',q'\in \N$ 
and we \lq \lq {\stampatello flip}\rq \rq \thinspace the divisors $\ell:=q/d$.) 
\par
{\stampatello the} $\underline{\hbox{\stampatello Ramanujan coefficients}}$ {\stampatello of our} $f:§\N \rightarrow \C$ {\stampatello are }
${\displaystyle 
R_{\ell}(f)\defineq \sum_{m\equiv 0(\!\!\bmod \ell)}{{(f\ast \mu)(m)}\over m} \enspace \forall \ell \in \N
}$
\par
({\stampatello well-defined, since} $g:=f\ast \mu$ {\stampatello has} $|\supporto(g)|<\infty$) {\stampatello to get} (with: $\supporto(g)\subset [1,Q]$ and $g\EssBdd 1$) 
$$
|g|\ll G\EssBdd 1 \Rightarrow 
R_{\ell}(g\ast \1)={1\over {\ell}}\sum_{q\le {Q\over {\ell}}}{{g(\ell q)}\over q}
\ll R_{\ell}(G\ast \1)
\EssBdd {1\over {\ell}}\sum_{q\le {{2N+h}\over {\ell}}}{1\over q}=R_{\ell}(\1 \ast \1):=R_{\ell}(d)\EssBdd {1\over {\ell}}. 
\leqno{(0)}
$$
\par
Here \thinspace $d(n)\defineq \sum_{q|n}1\EssBdd 1$ \thinspace is the {\stampatello divisor function}. We need the following Lemma. 

\bigskip
\bigskip
\bigskip

\par				
\noindent \centerline{\bf 2. An elementary Lemma.}
\smallskip
\par
\noindent
This inequality, saying \enspace $g\EssBdd 1$ $\Rightarrow$ $R_q(g\ast \1)\EssBdd R_q(\1 \ast \1)$\enspace ($\forall q$, here), is the core of the general philosophy underlying our {\stampatello theorem}: we use a kind of \lq \lq majorant\rq \rq, for (real) essentially bounded functions $f(n)$, represented by the divisor function, $d(n)$. However, for the time being, we use $|g|\ll G\EssBdd 1$ in our bounds. 
\par
Also, {\stampatello these coefficients allow us, then, to write} 
$$
\modSel_f(N,h)=\sum_{1<\ell \le Q}R_{\ell}^2(f)\starsum_{j\le {{\ell}\over 2}}\modFh^2\!\left({j\over {\ell}}\right) \sum_{x\sim N}\cos^2 {{2\pi xj}\over {\ell}} + 
$$
$$
          + 2\doublesum_{1<\ell,t \le Q}R_{\ell}(f)R_t(f)
              \doublestarsum_{{j\le {{\ell}\over 2} \enspace \thinspace \enspace r\le {t\over 2}}\atop {{j\over {\ell}}-{r\over t}>0}}
                 \modFh\left({j\over {\ell}}\right)\modFh\left({r\over t}\right)
                        \sum_{x\sim N}\cos{{2\pi xj}\over {\ell}}\cos{{2\pi xr}\over t} 
$$
\par
\noindent
(use above properties of $\chi_q$ and $(0)$ above), since for the {\stampatello Farey fractions} (they're reduced ones, this time in $[0,1/2]$) $j/\ell=r/t$ $\Rightarrow $ $j=r, \ell=t$ and we may exchange the couples whenever $j/\ell<r/t$. Hence, say, 
$$
\modSel_f(N,h)=\widetilde{D}_f(N,h) 
               +\doublesum_{1<\ell,t \le Q}R_{\ell}(f)R_t(f)
                 \doublestarsum_{{j\le {{\ell}\over 2} \enspace \thinspace \enspace r\le {t\over 2}}\atop {{j\over {\ell}}-{r\over t}>0}}
                   \modFh\left({j\over {\ell}}\right)\modFh\left({r\over t}\right)
                     \left(\sum_{x\sim N}\cos {2\pi \delta x}+\sum_{x\sim N}\cos {2\pi \sigma x}\right), 
\leqno{(1)}
$$
\par
\noindent
where, say, 
$$
\widetilde{D}_f(N,h)\defineq \sum_{1<\ell \le Q}R_{\ell}^2(f)\starsum_{j\le {{\ell}\over 2}}\modFh^2\!\left({j\over {\ell}}\right) \sum_{x\sim N}\cos^2 {{2\pi xj}\over {\ell}}\ge 0 
$$
\par
\noindent
is the {\stampatello diagonal} and, say, $\delta:=\left\Vert{j\over {\ell}}-{r\over t}\right\Vert={j\over {\ell}}-{r\over t}>0$, \enspace $\sigma:=\left\Vert{j\over {\ell}}+{r\over t}\right\Vert \in \left[0,{1\over 2}\right]$. Since $\modFh \,(j/q) \ge 0$ $\forall j\le q/2$, 
$$
\doublesum_{1<\ell,t \le Q}R_{\ell}(f)R_t(f)
                 \doublestarsum_{{j\le {{\ell}\over 2} \enspace \thinspace \enspace r\le {t\over 2}}\atop 
										{0<\delta:={j\over {\ell}}-{r\over t}\le {1\over A}}}
                   \modFh\left({j\over {\ell}}\right)\modFh\left({r\over t}\right)
                    \sum_{x\sim N}\cos {2\pi \delta x}
\EssBdd
$$
$$
\EssBdd
\doublesum_{1<\ell,t \le Q}R_{\ell}(G\ast \1)R_t(G\ast \1)
                 \doublestarsum_{{j\le {{\ell}\over 2} \enspace \thinspace \enspace r\le {t\over 2}}\atop 
										{0<\delta:={j\over {\ell}}-{r\over t}\le {1\over A}}}
                   \modFh\left({j\over {\ell}}\right)\modFh\left({r\over t}\right)
                    \sum_{x\sim N}\cos {2\pi \delta x}, 
$$
\par
\noindent
{\stampatello and} 
$$
\doublesum_{1<\ell,t \le Q}R_{\ell}(f)R_t(f)
                 \doublestarsum_{{j\le {{\ell}\over 2} \enspace \thinspace \enspace r\le {t\over 2}}\atop 
									{\delta>0,\sigma:=\left\Vert {j\over {\ell}}+{r\over t}\right\Vert\le {1\over A}}}
                   \modFh\left({j\over {\ell}}\right)\modFh\left({r\over t}\right)
                    \sum_{x\sim N}\cos {2\pi \sigma x}
\EssBdd
$$
$$
\EssBdd
\doublesum_{1<\ell,t \le Q}R_{\ell}(G\ast \1)R_t(G\ast \1)
                 \doublestarsum_{{j\le {{\ell}\over 2} \enspace \thinspace \enspace r\le {t\over 2}}\atop 
									{\delta>0,\sigma:=\left\Vert {j\over {\ell}}+{r\over t}\right\Vert\le {1\over A}}}
                   \modFh\left({j\over {\ell}}\right)\modFh\left({r\over t}\right)
                    \sum_{x\sim N}\cos {2\pi \sigma x}, 
$$
\par
\noindent
{\stampatello provided} $N=o(A)$, when $N\to \infty$. In fact, Taylor expansion of $\cos$ gives in both cases a positive $x-$sum.
\par
\noindent
Here we used estimates for Ramanujan coefficients coming from $(0)$ bounds. 
\medskip
\par
This \lq \lq {\stampatello majorant principle}\rq \rq \thinspace is not applicable to all of our $\modSel_f(N,h)$, as the terms for which $\delta$ or $\sigma$ are above $1/A$ ($\to 0$, say; better, $1/A=o(1/N)$ here) are troublesome: we don't know the sign of the $x-$sum.
\par
However, luckily enough, we are able bound their contribution to our integral using a very simple \lq \lq {\stampatello well-spaced}\rq \rq \thinspace argument, to be explicit the one used to prove Large Sieve type inequalities, applied to Farey fractions (as we have here, indeed). This has been done in Lemma 2 [C-S] (uses only Cauchy inequality); as in our present {\stampatello elementary} Lemma, following.

\medskip

\par				
We can state and show our
\smallskip
\par
\noindent {\stampatello Lemma.} {\it Let } $N,h\in \N$ {\it with } $h\to \infty$ {\it and } $h=o(N)$ {\it when } $N\to \infty$. {\stampatello Assume } $g:\N \rightarrow \C$ {\it with } $g(q)=0$ $\forall q>Q$, {\it where } $1\le Q\ll N$. {\it Set } $f:=g\ast \1$. {\it Choose } $A\in \R$, $A=A(N,h)\to \infty$ {\it when } $N\to \infty$. {\stampatello Then} 
$$
\doublesum_{1<\ell,t \le Q}R_{\ell}(f)R_t(f)
                 \doublestarsum_{{j\le {{\ell}\over 2} \enspace \thinspace \enspace r\le {t\over 2}}\atop {\delta:={j\over {\ell}}-{r\over t}>1/A}}
                   \modFh\left({j\over {\ell}}\right)\modFh\left({r\over t}\right)
                    \sum_{x\sim N}\cos {2\pi \delta x}\EssBdd Ah,
$$
$$
\doublesum_{1<\ell,t \le Q}R_{\ell}(f)R_t(f)
                 \doublestarsum_{{j\le {{\ell}\over 2} \enspace \thinspace \enspace r\le {t\over 2}}\atop 
									{\delta>0,\sigma:=\left\Vert{j\over {\ell}}+{r\over t}\right\Vert>1/A}}
                   \modFh\left({j\over {\ell}}\right)\modFh\left({r\over t}\right)
                    \sum_{x\sim N}\cos {2\pi \sigma x}\EssBdd Ah.
$$

\smallskip

\par
\noindent {\stampatello proof.}$\!$ The elementary calculation of the following exponential sum (compare [D, ch.25]) gives
$$
\alpha \not \in \Z \enspace \Rightarrow \enspace \sum_{x\sim N}e(\alpha x)\ll {1\over {\Vert \alpha \Vert}}, 
$$
\par
which, together with (recall $\modcjq:=\modFh\,(j/q)/q$, here use $(0)$ bounds, from $f\EssBdd 1$)
$$
\doublesum_{1<\ell,t \le Q}R_{\ell}(f)R_t(f)
                 \doublestarsum_{{j\le {{\ell}\over 2} \enspace \thinspace \enspace r\le {t\over 2}}\atop {\delta:={j\over {\ell}}-{r\over t}>1/A}}
                   \modFh\left({j\over {\ell}}\right)\modFh\left({r\over t}\right)
                    \sum_{x\sim N}\cos {2\pi \delta x}
\EssBdd
\doublesum_{1<\ell,t \le Q}\doublestarsum_{{|j|\le {{\ell}\over 2} \enspace \thinspace \enspace |r|\le {t\over 2}}\atop 
										{\left\Vert{j\over {\ell}}-{r\over t}\right\Vert>1/A}}
				{{|\widetilde{c}_{j,\ell}|\cdot |\widetilde{c}_{r,t}|}\over {\left\Vert{j\over {\ell}}-{r\over t}\right\Vert}}
$$
\par
and, changing sign to $r$, 
$$
\enspace 
\doublesum_{1<\ell,t \le Q}R_{\ell}(f)R_t(f)
                 \doublestarsum_{{j\le {{\ell}\over 2} \enspace \thinspace \enspace r\le {t\over 2}}\atop 
									{\delta>0,\sigma:=\left\Vert{j\over {\ell}}+{r\over t}\right\Vert>1/A}}
                   \modFh\left({j\over {\ell}}\right)\modFh\left({r\over t}\right)
                    \sum_{x\sim N}\cos {2\pi \sigma x}
\EssBdd
\doublesum_{1<\ell,t \le Q}\doublestarsum_{{|j|\le {{\ell}\over 2} \enspace \thinspace \enspace |r|\le {t\over 2}}\atop 
										{\left\Vert{j\over {\ell}}-{r\over t}\right\Vert>1/A}}
				{{|\widetilde{c}_{j,\ell}|\cdot |\widetilde{c}_{r,t}|}\over {\left\Vert{j\over {\ell}}-{r\over t}\right\Vert}},
$$
\par
{\stampatello give ${1\over A}$ well-spaced (Farey) fractions} (doesn't matter where, in $[-{1\over 2},{1\over 2}]$ or $[0,1]$ here); then, 
$$
\Sigma := 
\doublesum_{1<\ell,t \le Q}\doublestarsum_{{|j|\le {{\ell}\over 2} \enspace \thinspace \enspace |r|\le {t\over 2}}\atop 
										{\left\Vert{j\over {\ell}}-{r\over t}\right\Vert>1/A}}
				{{|\widetilde{c}_{j,\ell}|\cdot |\widetilde{c}_{r,t}|}\over {\left\Vert{j\over {\ell}}-{r\over t}\right\Vert}}
\ll
\sum_{1<\ell \le Q}\starsum_{j\le \ell}|\widetilde{c}_{j,\ell}|^2 
 \enspace \mathop{\sum_{1<t\le Q}\starsum_{r\le t}}_{\left\Vert {r\over t}-{j\over {\ell}}\right\Vert>1/A}{1\over {\left\Vert {r\over t}-{j\over {\ell}}\right\Vert}}, 
$$
\par
using Cauchy inequality (\& variables symmetry); number the ${\cal O}(Q^2)$ Farey fractions $\lambda_m:={j\over {\ell}}$, $\lambda_n:={r\over t}$, 
$$
n\neq m \Rightarrow \left\Vert \lambda_n-\lambda_m\right\Vert \ge {{|n-m|}\over A} 
$$
\par
which gives, recalling from the above 
$$
\starsum_{j\le \ell}|\widetilde{c}_{j,\ell}|^2 \ll \sum_{j\le \ell}|\widetilde{c}_{j,\ell}|^2 \ll \left\Vert {h\over {\ell}}\right\Vert \ll \min\left(1,{h\over {\ell}}\right), 
$$
\par
the required 
$$
\Sigma \EssBdd A\sum_{1<\ell \le Q}\starsum_{j\le \ell}|\widetilde{c}_{j,\ell}|^2 
\EssBdd Ah, 
$$
\par
since (in the sequel, let \enspace $L:=\log N$) 
$$
\mathop{\sum_{1<t\le Q}\starsum_{r\le t}}_{\left\Vert {r\over t}-{j\over {\ell}}\right\Vert>1/A}{1\over {\left\Vert {r\over t}-{j\over {\ell}}\right\Vert}}
=\sum_{{n\neq m}\atop {\left\Vert \lambda_n-\lambda_m\right\Vert>1/A}}{1\over {\left\Vert \lambda_n-\lambda_m\right\Vert}}
\ll A\sum_{n\neq m}{1\over {|n-m|}}\ll A\sum_{1\le k\ll Q^2}{1\over k}\ll AL 
\EssBdd A.\enspace \square
$$

\bigskip
\bigskip
\bigskip

\par				
\noindent \centerline{\bf 3. Completion of the Theorem proof. Proof of the Proposition. Remarks.}
\smallskip
\par
\noindent
We complete the proof, {\stampatello recalling the majorant principle}, see above, with \thinspace $f:=g\ast \1$, {\stampatello say}, \enspace $F:=G\ast \1$: 
$$
\widetilde{D}_f(N,h) 
               +\doublesum_{1<\ell,t \le Q}R_{\ell}(f)R_t(f)
                 \doublestarsum_{{j\le {{\ell}\over 2} \enspace \thinspace \enspace r\le {t\over 2}}\atop {0<\delta\le {1\over A}}}
                   \modFh\left({j\over {\ell}}\right)\modFh\left({r\over t}\right)
                     \sum_{x\sim N}\cos {2\pi \delta x} + 
\leqno{(2)}
$$
$$
               +\doublesum_{1<\ell,t \le Q}R_{\ell}(f)R_t(f)
                 \doublestarsum_{{j\le {{\ell}\over 2} \enspace \thinspace \enspace r\le {t\over 2}}\atop 
								{\delta>0,\sigma:=\left\Vert{j\over {\ell}}+{r\over t}\right\Vert \le {1\over A}}}
                   \modFh\left({j\over {\ell}}\right)\modFh\left({r\over t}\right)
                     \sum_{x\sim N}\cos {2\pi \sigma x}\EssBdd 
$$
$$
\EssBdd \widetilde{D}_F(N,h) 
               +\doublesum_{1<\ell,t \le Q}R_{\ell}(F)R_t(F)
                 \doublestarsum_{{j\le {{\ell}\over 2} \enspace \thinspace \enspace r\le {t\over 2}}\atop {0<\delta\le {1\over A}}}
                   \modFh\left({j\over {\ell}}\right)\modFh\left({r\over t}\right)
                     \sum_{x\sim N}\cos {2\pi \delta x} + 
$$
$$
               +\doublesum_{1<\ell,t \le Q}R_{\ell}(F)R_t(F)
                 \doublestarsum_{{j\le {{\ell}\over 2} \enspace \thinspace \enspace r\le {t\over 2}}\atop 
								{\delta>0,\sigma:=\left\Vert{j\over {\ell}}+{r\over t}\right\Vert \le {1\over A}}}
                   \modFh\left({j\over {\ell}}\right)\modFh\left({r\over t}\right)
                     \sum_{x\sim N}\cos {2\pi \sigma x}, 
$$
\par
{\stampatello say, for} ${1\over A}=o\left({1\over N}\right)$; {\stampatello since the (well-spaced Farey fractions) Lemma gives} plus $\EssBdd Ah$, {\stampatello both} 
\par
for $f$ and $F$ (in place of $f$), {\stampatello choose } $A=NL$ (or {\stampatello even } $A=N\log \log N$, here!) {\stampatello in order to get} 
$$
\modSel_{f}(N,h)\EssBdd \modSel_{G\ast \1}(N,h)+Nh.\enspace \square 
$$

\bigskip

\par
\noindent 
We prove, now, the Proposition. (Abbrev. $X(q):=${\stampatello the inverse} $Q^{-1}(q)$: $Q(X(q))=q$, $X(Q(x))=x$ $\forall q,x$.)
\smallskip
\par
\noindent {\stampatello proof.}$\!$ Instead of $(1)$ we have (recall from $(0)$ the $R_{\ell}(f)$ definition), {\stampatello with} \thinspace $Q:=Q(2N+h)$, that $\modSel_f(N,h)$ is 
$$
\enspace 
\sum_{x\sim N}\Big| \sum_{q\le Q(x+h)}g(q)\modchi_q(x)\Big|^2 
= \sum_{1<\ell \le Q} \doublesum_{d_1 \thinspace , \thinspace d_2 \thinspace \le {Q\over {\ell}}}{{g(\ell d_1)g(\ell d_2)}\over {d_1 d_2}}
					\starsum_{j\le {{\ell}\over 2}}\modFh^2\!\left( {j\over {\ell}}\right) 
						\sum_{{x\sim N}\atop {{x\ge X(\ell d_1) - h}\atop {x\ge X(\ell d_2) - h}}}\cos^2 {{2\pi xj}\over {\ell}} + 
$$
$$
\enspace 
               +\doublesum_{1<\ell,t \le Q}\sum_{d\le {Q\over {\ell}}}\sum_{q\le {Q\over t}}{{g(\ell d)g(tq)}\over {dq}}
                 \doublestarsum_{{j\le {{\ell}\over 2} \enspace \thinspace \enspace r\le {t\over 2}}\atop {\delta:={j\over {\ell}}-{r\over t}>0}}
                   \modFh\left({j\over {\ell}}\right)\modFh\left({r\over t}\right)
                     \left(\sum_{{x\sim N}\atop {{x\ge X(\ell d)-h}\atop {x\ge X(tq)-h}}}\cos {2\pi \delta x} 
				+ \sum_{{x\sim N}\atop {{x\ge X(\ell d)-h}\atop {x\ge X(tq)-h}}}\cos {2\pi \sigma x}\right) 
$$
\par
and \enspace $(i)$ the \lq \lq {\stampatello majorant part}\rq \rq \thinspace ({\stampatello i.e., diagonal \& nearby, $(2)$ above}) is {\stampatello unchanged}, since {\stampatello Taylor} 
\par
{\stampatello expansion} still gives {\stampatello non-negative $x-$sums}, while \enspace $(ii)$ the \lq \lq {\stampatello well-spaced part}\rq \rq, \thinspace i.e., {\stampatello the Lemma}, 
\par
{\stampatello still holds}, since in our {\stampatello Lemma} the initial estimate (of {\stampatello exponential sums, from} [D], quoted) {\stampatello does}
\par
{\stampatello not depend on the $x-$interval} ({\stampatello though maybe non-optimal}, say the length is $o(1/\Vert \alpha\Vert)$, esp.).\enspace $\square$ 

\bigskip

\par
\noindent
We give some remarks, to have an idea of the (possible) applications of the Theorem (and/or the Proposition). 
\smallskip
\par
First of all, we may think to substitute $G$ with a constant (say, $G=\1$), but then the resulting modified Selberg integral of the divisor function does not have a non-trivial estimate, since we now force the expected value $M_f(x,h)/h$ of $f$ in the short interval to be of the kind $\sum_{d\le x}{1\over d}$, which is different from the actual one (roughly $\log x + 2\gamma$, $\gamma:=${\stampatello Euler-Mascheroni constant}, instead of the present $\log x + \gamma$); this results in a \lq \lq trivial\rq \rq \thinspace modified Selberg integral. In fact, if we look (say) for example, at [C-S] Theorem 2 (even if about Selberg integral of the divisor function, not the modified one!), the expected value is calculated (in full agreement with classic residue-calculated terms, see [C4], for example) starting from $g=\1$, but after flipping the divisors (so, the range is not up to about $x$, but cut  about $\sqrt x$ !). 
\par				
Second, the search for good (i.e., with non-trivial modified Selberg integral $\modSel_{G\ast \1}(N,h)$, say we gain small powers on the rough ${\cal O}(Nh^2)-$bound) majorants $G$, i.e. good Wintner majorants $F:=G\ast \1$ for the (real, essentially bounded) $f:=g\ast \1$ (i.e., $|g|\ll G$), is not a trivial question !
\par
Third, it's not yet completely clear that this can (and how) give a \lq \lq smoothing\rq \rq \thinspace of the arithmetic behind the function $f:=g\ast \1$ (substantially and morally, we should like to bound $g$ with constants, but see the first consideration above !); compare (but there we make further hypotheses on the function $f$) the appearance of the idea of majorant principles in our paper [C5] (where the attention is on non-negative exponential sums in the long range, not in the short, that's for free here !).
\par
Last but not least, why did we study the modified Selberg integral and not simply the Selberg integral with this approach ? Simply because it fails, due to Fourier coefficients of non-constant sign coming from the short interval (see, instead, the coefficients $\modFh \ge 0$ above) ! What about the symmetry integral, then ...

\bigskip
\bigskip
\bigskip

\par
\noindent
\centerline{\bf References}
\smallskip
\item{\bf [C1]} \thinspace Coppola, G.\thinspace - \thinspace {\sl On the Correlations, Selberg integral and symmetry of sieve functions in short intervals} \thinspace - \thinspace http://arxiv.org/abs/0709.3648v3 (to appear on: Journal of Combinatorics and Number Theory)
\smallskip
\item{\bf [C2]} \thinspace Coppola, G.\thinspace - \thinspace {\sl On the Correlations, Selberg integral and symmetry of sieve functions in short intervals, II} \thinspace - \thinspace Int. J. Pure Appl. Math. {\bf 58.3}(2010), 281--298.
\smallskip
\item{\bf [C3]} \thinspace Coppola, G.\thinspace - \thinspace {\sl On the Correlations, Selberg integral and symmetry of sieve functions in short intervals, III} \thinspace - \thinspace http://arxiv.org/abs/1003.0302v1
\smallskip
\item{\bf [C4]} \thinspace Coppola, G.\thinspace - \thinspace {\sl On the Selberg integral of the $k-$divisor function and the $2k-$th moment of the Riemann zeta-function} \thinspace - \thinspace http://arxiv.org/abs/0907.5561v1 - to appear on Publ. Inst. Math., Nouv. Sér.
\smallskip
\item{\bf [C5]} \thinspace Coppola, G.\thinspace - \thinspace {\sl On the symmetry of arithmetical functions in almost all short intervals, V} \thinspace - \thinspace (electronic) http://arxiv.org/abs/0901.4738v2 
\smallskip
\item{\bf [C-S]} Coppola, G. and Salerno, S.\thinspace - \thinspace {\sl On the symmetry of the divisor function in almost all short intervals} \thinspace - \thinspace Acta Arith. {\bf 113} (2004), {\bf no.2}, 189--201. $\underline{\tt MR\enspace 2005a\!:\!11144}$
\smallskip
\item{\bf [D]} \thinspace Davenport, H.\thinspace - \thinspace {\sl Multiplicative Number Theory} \thinspace - \thinspace Third Edition, GTM 74, Springer, New York, 2000. $\underline{{\tt MR\enspace 2001f\!:\!11001}}$
\smallskip
\item{\bf [T]} \thinspace Tenenbaum, G.\thinspace - \thinspace {\sl Introduction to Analytic and Probabilistic Number Theory} \thinspace - \thinspace Cambridge Studies in Advanced Mathematics, {\bf 46}, Cambridge University Press, 1995. $\underline{\tt MR\enspace 97e\!:\!11005b}$
\smallskip
\item{\bf [V]} \thinspace Vinogradov, I.M.\thinspace - \thinspace {\sl The Method of Trigonometrical Sums in the Theory of Numbers} - Interscience Publishers LTD, London, 1954. $\underline{{\tt MR \enspace 15,941b}}$

\medskip

\leftline{\tt Dr.Giovanni Coppola}
\leftline{\tt DIIMA - Universit\`a degli Studi di Salerno}
\leftline{\tt 84084 Fisciano (SA) - ITALY}
\leftline{\tt e-mail : gcoppola@diima.unisa.it}
\leftline{\tt e-page : www.giovannicoppola.name}

\bye